\documentclass[10pt]{amsart}
\usepackage{latexsym, amsmath, amssymb}
\usepackage{epsfig,graphics}

\newcommand{\newsection}[1]
{\section{#1}\setcounter{theorem}{0} \setcounter{equation}{0}
\par\noindent}

\newtheorem{theorem}{Theorem}

\newcommand{\cd}{\, \cdot\, }

\newcommand{\R}{{\mathbb R}}

\newcommand{\ang}{{\not\negmedspace\nabla}}

\newcommand{\la}{\langle}
\newcommand{\ra}{\rangle}

\newcommand{\M}{{\mathcal{M}}}

\newcommand{\weight}{{\Bigl(1-\Bigl(\frac{r_s}{r}\Bigr)^{d+1}\Bigr)}}
\renewcommand{\S}{{\mathbb{S}}}

\begin{document}

\title[Localized energy estimates on Schwarzschild space-times]
{
Localized energy estimates for wave equations on high dimensional
Schwarzschild space-times
}

\thanks{
The second author was supported in part by the NSF through grant DMS0800678
}

\author{Parul Laul}
\author{Jason Metcalfe}

\address{Department of Mathematics, University of North Carolina,
  Chapel Hill, NC 27599-3250}

\begin{abstract}
The localized energy estimate for the wave equation is known to be a
fairly robust measure of dispersion.  Recent analogs on
the $(1+3)$-dimensional Schwarzschild space-time have played a key
role in a number of subsequent results, including a proof of Price's
law.  In this article, we explore similar localized energy estimates
for wave equations on $(1+n)$-dimensional hyperspherical Schwarzschild
space-times.
\end{abstract}

\maketitle


\newsection{Introduction}
One of the more robust measures of dispersion for the wave equation
is the so called localized energy estimates.  These estimates have
played a role in understanding scattering theory, as means of summing
local in time Strichartz estimates to obtain global in time
Strichartz estimates on a compact set, as an essential tool for
proving long time existence of quasilinear wave equations in exterior
domains, and as means to handle errors in certain parametrix
constructions.  

Estimates of this form originated in \cite{M} where it was shown that
solutions to the constant coefficient, homogeneous wave equation
\[(\partial_t^2-\Delta)u=0,\quad u(0)=u_0,\quad \partial_t u(0)=u_1\]
on $\R_+\times \R^n$, $n\ge 3$ satisfy
\[\int_0^T \int_{\R^n} \frac{1}{|x|}|\ang u|^2(t,x)\,dx\,dt\lesssim
\|\nabla u_0\|_2^2 + \|u_1\|_2^2\]
where $\ang$ denotes the angular derivatives.  Though not the original
method of proof, one can obtain this estimate by multiplying $\Box u$
by $\partial_r u +\frac{n-1}{2}\frac{u}{r}$, integrating over a
space-time slab, and integrating by parts.  One can also obtain
control on $\partial_r$ and $\partial_t$ in a fixed dyadic annulus.
Attempting to sum over these dyadic annuli comes at the cost of a
logarithmic blow up in time.  One may otherwise introduce an
additional component of the weight that permits summability.  In this
case, it can be shown that for $u$ as above and $n\ge 4$
\begin{equation}\label{kss}
\|\la x\ra^{-1/2-} u'\|_{L^2_{t,x}([0,T]\times\R^n)} + \|\la
x\ra^{-3/2} u\|_{L^2_{t,x}([0,T]\times \R^n)} \lesssim \|\nabla u_0\|_2
+ \|u_1\|_2.
\end{equation}
Here $\la x\ra = \sqrt{1+|x|^2}$ denotes the Japanese bracket and
$u'=(\partial_t u,\nabla_xu)$ represents the space-time gradient.  An
estimate akin to \eqref{kss} also holds for $n=3$ but in this case the
weight $\la x\ra^{-3/2}$ in the second term in the left side must be
replaced by $\la x\ra^{-3/2-}$.  To prove \eqref{kss}, one multiplies
the equation instead by $f(r)\partial_r u +
\frac{n-1}{2}\frac{f(r)}{r}u$ where $f(r)=\frac{r}{r+2^j}$, $j\ge 0$
and integrates by parts.  For fixed $j$, this yields an estimate for
the left side over $|x|\approx 2^j$ with weight $\la x\ra^{-1/2}$ in
the first term.  Introducing the additional weight and summing over
$j$ produces \eqref{kss}.  See, e.g., \cite{HY}, \cite{KSS}, \cite{KPV},
\cite{Met}, \cite{SmSo}, \cite{St}, \cite{Strauss}.  Analogous
estimates have been shown for small, asymptotically flat, possibly time-dependent
perturbations of the d'Alembertian in \cite{Alinhac}, \cite{MS},
\cite{MT, MT2} as well as for time-independent, asymptotically flat, nontrapping
perturbations in e.g. \cite{BH}, \cite{Burq}, \cite{SW}.

A particular case of interest which does not fall into these latter
categories is the wave equation on (1+3)-dimensional Schwarzschild space-times.  While
asymptotically flat and time independent, this metric is not, however,
nontrapping.  There is trapping which occurs on the so called photon
sphere which necessitates a loss in the estimates as compared to those
for the Minkowski wave equation.  Despite this complication, a number
of proofs of estimates akin to \eqref{kss} exist for the wave equation
on Schwarzschild space-times.  See \cite{BS1}-\cite{BS5}, \cite{DR, DR2},
\cite{MMTT}.  In \cite{TT}, \cite{DR5}, and \cite{AB}, these estimates have been extended to Kerr
space-times with small angular momentums.  Here an additional
difficulty is encountered as \cite{Alinhac2} shows that the estimates
will not follow from an analog of the argument above with any first order
differential multiplier.

The goal of this article is to extend the known localized energy
estimates on $(1+3)$-dimensional Schwarzschild space-times to
$(1+n)$-dimensional Schwarzschild space-times for $n\ge 3$.  There are
multiple notions of the Schwarzschild space-time for $n\ge 4$, and we
restrict our attention to the hyperspherical case.  For a derivation
and discussion of such black hole space times, see
e.g. \cite{T}, \cite{MP}.  An alternative notion of higher
dimensional Schwarzschild space-times is discussed e.g. in \cite{GL}.
We shall not explore these hypercylindrical Schwarzschild manifolds or
other notions of higher dimensional black holes here.

The exterior of a $(1+n)$-dimensional hyperspherical Schwarzschild
black hole $(\M,g)$ is described by the manifold $\M=\R\times
(r_s,\infty)\times \S^{d+2}$ and the line element
\[ds^2 = -\weight dt^2 + \weight^{-1}dr^2 + r^2d\omega^2.\]
Here $d=n-3$, $r_s$ denotes the Schwarzschild radius (i.e. $r=r_s$ is
the event horizon), and $d\omega$ is
the surface measure on the sphere $\S^{d+2}=\S^{n-1}$.  Here, the
d'Alembertian is given by
\begin{align*}
  \Box_g \phi &=\nabla^\alpha \partial_\alpha \phi\\
&= -\weight^{-1}\partial_t^2 \phi + r^{-(d+2)}\partial_r
  \Bigl[r^{d+2}\weight\partial_r \phi\Bigr] + \ang\cdot \ang \phi.
\end{align*}
The Killing vector field $\partial_t$ yields the
conserved energy
\begin{multline*}
E[\phi](t)=\int_{\S^{d+2}}\int_{r\ge r_s}
\Bigl[\weight^{-1}(\partial_t\phi)^2(t,r,\omega) \\+ \weight
  (\partial_r\phi)^2(t,r,\omega) + |\ang \phi|^2(t,r,\omega)\Bigr]r^{d+2}\,dr\,d\omega.
\end{multline*}
That is, when $\Box_g\phi = 0$, $E[\phi](t)=E[\phi](0)$ for all $t$.

We now define our localized energy norm.  To this end, we set
\begin{multline}\label{LEnorm}
  \|\phi\|_{LE}^2 = \int_0^\infty \int_{\S^{d+2}}\int_{r\ge r_s} \Bigl[c_r(r) \weight
  (\partial_r\phi)^2(t,r,\omega)
+ c_\omega(r) |\ang\phi|^2(t,r,\omega) \\+
c_0(r)\phi^2(t,r,\omega)\Bigr] r^{d+2}\,dr\,d\omega\,dt
\end{multline}
where 
\[
c_r=\frac{1}{r^{d+3}\Bigl(1-\log\Bigl(\frac{r-r_s}{r}\Bigr)\Bigr)^2},\quad
c_\omega=\frac{1}{r} \Bigl(\frac{r-r_{ps}}{r}\Bigr)^2,\]
\[c_0=\Bigl(\frac{r-r_s}{r}\Bigr)^{-1}\frac{1}{r^3\Bigl(1-\log\Bigl(\frac{r-r_s}{r}\Bigr)\Bigr)^4}. \]
Here $r_{ps}=\Bigl(\frac{d+3}{2}\Bigr)^{\frac{1}{d+1}}r_s$, which is
the location of the photon sphere.
The main result of this article then states that a localized energy
estimate holds on $\M$.
\begin{theorem}\label{thm1.1}\footnotemark
  Let $\phi$ solve the homogeneous wave equation $\Box_g\phi=0$ on a
  $(1+n)$-dimensional hyperspherical Schwarzschild manifold with $n\ge
  4$.  Then we have
\begin{equation}\label{LEestimate}
\sup_{t\ge 0}E[\phi](t)+\|\phi\|_{LE}^2\lesssim E[\phi](0).\end{equation}
\end{theorem}

\footnotetext{Upon completion of this work, the authors learned that a
version of Theorem \ref{thm1.1} was also obtained independently in the
forthcoming \cite{Schlue}.}

We note that by modifying the lower order portion of the multiplier
which appears in the next section that it is straightforward to obtain an
estimate on $\partial_t\phi$ as well.  Moreover, one can obtain an
analogous estimate for the inhomogeneous wave equation by placing the
forcing term in an appropriate dual norm.  The necessary decay of the
coefficient $c_r$ at $\infty$ can also be significantly improved.  It is relatively
simple to carry out these modifications, but for the sake of
clarity, we omit the details.

As in the $(1+3)$-dimensional case, the higher dimensional
hyperspherical Schwarzschild space-times have a photon sphere where
trapping occurs.  Rays initially located on this sphere and moving
initially tangent to it will stay on the surface for all times.  Such
trapping is an obstacle to many types of dispersive estimates, and the
vanishing in the coefficient $c_\omega$ at the trapped set is a loss,
when compared to the nontrapping Minkowski setting, which results from
this trapping.

The proof follows by constructing an appropriate differential
multiplier.  This is carried out in the next section, and the
construction most closely resembles that which appears in
\cite{MMTT}.  Some significant technicalities needed to be
resolved, however, in order to find a construction which works in all
dimensions $d\ge 1$.

The localized energy estimates on $(1+3)$-dimensional Schwarzschild
manifolds have played a key role in a number of subsequent results.
See, e.g., \cite{Blue}, \cite{BS1}-\cite{BS5},
\cite{BSt},\cite{DR4}-\cite{DR3}, \cite{Luk}, \cite{MMTT},
\cite{Tataru}, \cite{TT}.  Theorem \ref{thm1.1} permits possible
generalizations of these studies to higher dimensions, and a portion
of these will be explored in the first author's upcoming doctoral
dissertation and other subsequent works.

\newsection{Construction of the multiplier}
Associated to $\Box_g$ is the energy-momentum tensor
\[Q_{\alpha\beta}[\phi]=\partial_\alpha\phi\partial_\beta\phi -
\frac{1}{2}g_{\alpha\beta}\partial^\gamma \phi \partial_\gamma \phi\]
whose most important property is the following divergence condition
\[\nabla^\alpha Q_{\alpha\beta}[\phi] = \partial_\beta\phi \Box_g
\phi.\]
Contracting $Q_{\alpha\beta}$ with the radial vector field
$X=f(r)\weight\partial_r$, we form the momentum density
$P_\alpha[\phi,X]=Q_{\alpha\beta}[\phi]X^\beta$.
Calculating the divergence of this quantity we have 
\begin{multline*}
  \nabla^\alpha P_\alpha[\phi,X] = \Box_g \phi X\phi +
  f'(r)\weight^2(\partial_r\phi)^2
+
\Bigl(\frac{r^{d+1}-\frac{d+3}{2}r_s^{d+1}}{r^{d+1}}\Bigr)\frac{f(r)}{r}|\ang
\phi|^2 
\\-\frac{1}{2}\Bigl[\weight r^{-(d+2)}\partial_r
  (r^{d+2}f(r))\Bigr]\partial^\gamma\phi\partial_\gamma \phi.
\end{multline*}
The last term involving the Lagrangian is not signed.  In order to
eliminate it, we utilize a lower order term in the multiplier.  To
this end, we modify the momentum density 
\begin{multline*}
\tilde{P}_\alpha[\phi,X] = P_\alpha[\phi,X] + \frac{1}{2}\Bigl[\weight r^{-(d+2)}\partial_r
  (r^{d+2}f(r))\Bigr]\phi\partial_\alpha \phi \\-
\frac{1}{4}\partial_\alpha\Bigl[\weight r^{-(d+2)}\partial_r
  (r^{d+2}f(r))\Bigr] \phi^2,
\end{multline*}
and recompute the divergence
\begin{multline}\label{div}
  \nabla^\alpha \tilde{P}_\alpha[\phi,X]=\Box_g\phi\Bigl[X\phi +
    \frac{1}{2}\Bigl\{\weight
    r^{-(d+2)}\partial_r(f(r)r^{d+2})\Bigr\}\Bigr]
\\+\weight^2 f'(r)(\partial_r\phi)^2 + \Bigl(\frac{r^{d+1} 
-r_{ps}^{d+1}}{r^{d+1}}\Bigr)\frac{f(r)}{r}|\ang \phi|^2 
\\-\frac{1}{4}\nabla^\alpha\partial_\alpha\Bigl[\weight
  r^{-(d+2)}\partial_r(f(r)r^{d+2})\Bigr] \phi^2.
\end{multline}

Ideally one would choose $f(r)$ to be smooth, bounded, and so that the
coefficients in the last three terms on the right are all nonnegative.
Indeed, this is precisely what is done in, e.g., \cite{St} and
\cite{MS} when $r_s=0$, i.e. in Minkowski space-time.  The localized
energy estimate then follows from an application of the divergence
theorem as well as an application of a Hardy inequality and the energy
estimate in order to handle the time boundary terms.
Unfortunately,
in the current setting it does not appear possible to construct such
an $f$.  We shall instead construct $f$ so that the last term in the right can be
bounded below by a positive quantity minus a fractional multiple of
the second term in the right.

We let $g(r)=\frac{r^{d+2}-r_{ps}^{d+2}}{r^{d+2}}$ and
$h(r)=\ln\Bigl(\frac{r^{d+1}-r_s^{d+1}}{\frac{d+1}{2}r_s^{d+1}}\Bigr)$.
The multiplier will be defined piecewise, and
it will be convenient to parametrize in terms of the values of
$h(r)$.  To this end, we fix the notation $r_\theta$ to denote the
value of $r$ so that $h(r_\theta)=\theta$.  Note, e.g., that
$r_{-\infty}=r_s$ and $r_{ps}=r_0$.  More explicitly, $r_\theta^{d+1}
=r_s^{d+1}\Bigl(\frac{d+1}{2}e^\theta + 1\Bigr)$.   An approximate multiplier is
\[g(r)+\frac{d+2}{d+3}\frac{r_{ps}r_s^{d+1}}{r^{d+2}}h(r).\]
We must, however, smooth out the logarithmic blow up at $r=r_s$.  We
must also smooth out $h(r)$ near $\infty$ in order to prevent a term
in $f'(r)$ which has an unfavorable sign.

In order to accomplish this, 
set
\[a(x)=
\begin{cases}
  -\frac{1}{\varepsilon} \frac{\varepsilon x + 1}{\delta(\varepsilon x+1)-1} -
  \frac{1}{\varepsilon},&x\le -\frac{1}{\varepsilon}\\
x, &-\frac{1}{\varepsilon}\le x \le 0\\
x-\frac{2}{3\alpha^2} x^3 +\frac{1}{5\alpha^4}x^5, &0\le x\le
\alpha\\
\frac{8\alpha}{15},&x\ge \alpha.
\end{cases}
\]
Here $\alpha = 5-\delta_0$ for some $0<\delta_0\ll 1$.  Then,
set
\[f(r)=g(r)+\frac{d+2}{d+3}\frac{r_{ps}r_s^{d+1}}{r^{d+2}}a(h(r)).\]
We notice that $f$ is $C^2$ with the exception of a jump in the second
derivative at $r_{-1/\varepsilon}$.

Using that $\Box_g\phi=0$, integrating \eqref{div} over $[0,T]\times
(r_s,\infty)\times \S^{d+2}$, and applying the divergence theorem, we have
\begin{multline}\label{base}
  -\int\int f(r)\partial_t\phi \partial_r\phi\,
  r^{d+2}\,dr\,d\omega\Bigl|_0^T - \frac{1}{2}\int\int
  \frac{1}{r^{d+2}}\partial_r(f(r)r^{d+2}) \phi \partial_t\phi\,
  r^{d+2}\,dr\,d\omega\Bigl|_0^T
\\=\int_0^T\int\int \Bigl\{ \weight^2 f'(r)(\partial_r\phi)^2 + 
 \Bigl(\frac{r^{d+1} 
-r_{ps}^{d+1}}{r^{d+1}}\Bigr)\frac{f(r)}{r}|\ang \phi|^2 
\\+l(f)\phi^2\Bigr\}r^{d+2}\,dr\,d\omega\,dt
\\+\frac{1}{4}r_{-1/\varepsilon}^{d+2}\Bigl(1-\Bigl(\frac{r_s}{r_{-1/\varepsilon}}\Bigr)^{d+1}\Bigr)^2
(f''(r_{-1/\varepsilon}^-)-f''(r_{-1/\varepsilon}^+))\int_0^T\int \phi^2|_{r=r_{-1/\varepsilon}}\,d\omega\,dt
\end{multline}
where
\[
l(f)=-\frac{1}{4}r^{-(d+2)}\partial_r\Bigl[\weight r^{d+2}\partial_r\Bigl\{\weight
  r^{-(d+2)}\partial_r(f(r)r^{d+2})\Bigr\}\Bigr].
\]

We first make a note about the boundary terms at
$r_{-1/\varepsilon}$.  Indeed, an elementary calculation shows
\begin{align*}
f''(r_{-1/\varepsilon}^-)-f''(r_{-1/\varepsilon}^+) &=\frac{d+2}{d+3}\frac{r_{ps}r_s^{d+1}}{r_{-1/\varepsilon}^{d+2}}
a''(-1/\varepsilon^-)(h'(r_{-1/\varepsilon}))^2 \\&= 2\delta\varepsilon
\frac{d+2}{d+3}\frac{r_{ps}r_s^{d+1}}{r_{-1/\varepsilon}^{d+2}}
(h'(r_{-1/\varepsilon}))^2\approx \delta \varepsilon e^{2/\varepsilon}.\end{align*}
Of particular interest is that the coefficient of the resulting boundary term at
$r_{-1/\varepsilon}$ is $O(\varepsilon\delta)$.

We now proceed to showing that the sum of the first three terms in the
right of \eqref{base} produce a positive contribution.  To do so, we will examine $f$
on a case by case basis and show 
\begin{itemize}
\item $f'(r)>0$ for all $r>r_s$
\item $f(r)<0$ for $r<r_{ps}$ and $f(r)>0$ for $r>r_{ps}$.
\item $\int l(f)\phi^2 \,r^{d+2}\,dr\,d\omega\,dt$ is bounded
  below by a positive term minus a fractional multiple of the
  $(\partial_r\phi)^2$ term and a $r_{-1/\varepsilon}$ boundary
  term.
\end{itemize}
By absorbing these latter pieces into those previously shown to
positively contribute, the estimate shall nearly be in hand.  It will
only remain to examine the time boundary terms, and in particular,
establish the Hardy inequality which will permit a direct application
of conservation of energy.

\noindent{\bf\em Case 1:} $r_s\le r\le r_{-1/\varepsilon}$.

The multiplier in this region is constructed to smooth out the
logarithmic blow up at the event horizon.  This is the only case in
which we shall not be able to just show that $l(f)\ge 0$.

Noting that $a(x)\le -1/\varepsilon$ for $x\le -1/\varepsilon$, we
immediately see that $f(r)<0$ on this range.  We also compute
\begin{multline}\label{df}
f'(r)=\frac{(d+2)r_{ps}^{d+2}}{r^{d+3}}-\frac{(d+2)^2}{d+3}\frac{r_{ps}r_s^{d+1}}{r^{d+3}}a(h(r))
\\+ \frac{d+2}{d+3}\frac{r_{ps}r_s^{d+1}}{r^{d+2}}\frac{1}{(\delta\varepsilon
  h(r)+\delta-1)^2}\frac{(d+1)r^d}{r^{d+1}-r_s^{d+1}}.
\end{multline}
Each of these summands is nonnegative on the given range, yielding the
desired sign for $f'(r)$.

It remains to show an appropriate lower bound for the $l(f)\phi^2$
term of \eqref{base}.  Calculating $l(f)$ we find
\begin{multline}\label{l_f}
 l(g(r))+l\Bigl(\frac{d+2}{d+3}\frac{r_{ps}r_s^{d+1}}{r^{d+2}}a(h(r))\Bigr)
 = \frac{d+2}{4r^{2d+5}}\Bigl(d
 r^{2d+2}+(d+3)r_s^{d+1}r^{d+1}-(d+2)^2r_s^{2d+2}\Bigr)
\\+\frac{(d+1)(d+2)}{2}\frac{r_{ps}r_s^{d+1}}{r^{2d+6}}(r_{ps}^{d+1}-r^{d+1})a'(h(r))
\\+\frac{(d+1)^2(d+2)(d+5)}{4(d+3)}\frac{r_{ps}r_s^{d+1}}{r^{d+5}}a''(h(r))
\\-\frac{(d+1)^3(d+2)}{4(d+3)}\frac{r_{ps}r_s^{d+1}}{r^4} \frac{1}{r^{d+1}-r_s^{d+1}}a'''(h(r)).
\end{multline}
As $a'(h(r))=(\delta\varepsilon h(r)+\delta -1)^{-2}$ and $r<r_{ps}$ in this regime,
the second term in the right has the desired sign.  The third term in
the right of \eqref{l_f} also has the desired sign as
$a''(h(r))=-2\delta\varepsilon/(\varepsilon \delta h(r) +\delta -1)^3$
and $h(r)\le -1/\varepsilon$ here.

The key step is to control the contribution of the last term in the
right of \eqref{l_f}.  We shall abbreviate $R(r)=\delta\varepsilon h(r)+\delta-1$.
By the Fundamental Theorem of Calculus, we observe
\[\int_{r_s}^{r_{-1/\varepsilon}} \partial_r\Bigl(\frac{2\delta\varepsilon
  (d+1)^2(d+2)}{d+3}\frac{r_{ps}r_s^{d+1}}{r^2(R(r))^3}\phi^2\Bigr)\,dr
=-\frac{2\delta\varepsilon
  (d+1)^2(d+2)}{d+3}\frac{r_{ps}r_s^{d+1}}{r_{-1/\varepsilon}^2}
\phi(r_{-1/\varepsilon})^2.\]
Evaluating the derivative in the left side yields
\begin{multline}\label{ineq1}
\int_{r_s}^{r_{-1/\varepsilon}} 
\frac{6\delta^2\varepsilon^2(d+1)^3(d+2)}{d+3}\frac{r_{ps}r_s^{d+1}}{r^2(R(r))^4}\frac{r^d}{r^{d+1}-r_s^{d+1}}\phi^2
\,dr
\\=-\int_{r_s}^{r_{-1/\varepsilon}} \frac{4\delta\varepsilon(d+1)^2(d+2)}{
  d+3}\frac{r_{ps}r_s^{d+1}}{r^3(R(r))^3}\phi^2\,dr
\\+\int_{r_s}^{r_{-1/\varepsilon}}
  \frac{4\delta\varepsilon(d+1)^2(d+2)}{d+3}\frac{r_{ps}r_s^{d+1}}{r^2(R(r))^3}\phi\partial_r \phi\,dr
\\+\frac{2\delta \varepsilon
  (d+1)^2(d+2)}{d+3}\frac{r_{ps}r_s^{d+1}}{r_{-1/\varepsilon}^2}
\phi(r_{-1/\varepsilon})^2.\end{multline}
To the second term in the right, we apply the Schwarz inequality to
obtain
\begin{multline}\label{ineq2}
  \int_{r_s}^{r_{-1/\varepsilon}}
  \frac{4\delta\varepsilon(d+1)^2(d+2)}{d+3}\frac{r_{ps}r_s^{d+1}}{r^2(R(r))^3}\phi\partial_r
  \phi\,dr
\\\le \int_{r_s}^{r_{-1/\varepsilon}} \frac{3\delta^2\varepsilon^2(d+2)(d+1)^3}{
  d+3} \frac{r_{ps}r_s^{d+1}r^{d-2}}{(R(r))^4(r^{d+1}-r_s^{d+1})}\phi^2\,dr
\\+\frac{4}{3}\int_{r_s}^{r_{-1/\varepsilon}} \frac{(d+2)(d+1)}{
  d+3} \frac{r_{ps}r_s^{d+1}}{r^{d+2}}\frac{r^{d+1}-r_s^{d+1}}{(R(r))^2}(\partial_r \phi)^2\,dr.
\end{multline}
Plugging \eqref{ineq2} into \eqref{ineq1} yields
\begin{multline}\label{combined}
  \int_{r_s}^{r_{-1/\varepsilon}} 
\frac{3\delta^2\varepsilon^2(d+1)^3(d+2)}{d+3}
\frac{r_{ps}r_s^{d+1}}{r^2(R(r))^4}\frac{r^d}{r^{d+1}-r_s^{d+1}}\phi^2
\,dr
\\\le -\int_{r_s}^{r_{-1/\varepsilon}} \frac{4\delta\varepsilon(d+1)^2(d+2)}{
  d+3}\frac{r_{ps}r_s^{d+1}}{r^3(R(r))^3}\phi^2\,dr
\\+\frac{4}{3}\int_{r_s}^{r_{-1/\varepsilon}} \frac{(d+2)(d+1)}{
  d+3} \frac{r_{ps}r_s^{d+1}}{r^{d+2}}\frac{r^{d+1}-r_s^{d+1}}{(R(r))^2}(\partial_r \phi)^2\,dr
\\+\frac{2\delta \varepsilon
  (d+1)^2(d+2)}{d+3}\frac{r_{ps}r_s^{d+1}}{r_{-1/\varepsilon}^2}
\phi(r_{-1/\varepsilon})^2.
\end{multline}
This shows that
\begin{multline}\label{rewritten}
  \Bigl(\frac{1}{4}+\frac{1}{48}\Bigr)\int_{\{r\in
    [r_s,r_{-1/\varepsilon}]\}}
\frac{(d+1)^3(d+2)}{d+3}\frac{r_{ps}r_s^{d+1}}{r^4}\frac{1}{r^{d+1}-r_s^{d+1}}a'''(h(r))\phi^2
\,r^{d+2}\,dr\,d\omega\,dt
\\\le \frac{13}{12}\int_{\{r\in [r_s,r_{-1/\varepsilon}]\}} \frac{(d+1)^2(d+2)}{
  (d+3)}\frac{r_{ps}r_s^{d+1}}{r^{d+5}}a''(h(r))\phi^2\,r^{d+2}\,dr\,d\omega\,dt
\\+\frac{13}{18}\int_{\{r\in[r_s,r_{-1/\varepsilon}]\}} \frac{(d+2)(d+1)}{
  (d+3)}
\frac{r_{ps}r_s^{d+1}}{r^{d+3}}\Bigl(1-\frac{r_s^{d+1}}{r^{d+1}}\Bigr)\frac{1}{(
  R(r))^2}(\partial_r \phi)^2\,r^{d+2}\,dr\,d\omega\,dt
\\+\frac{13}{12}\frac{\delta \varepsilon
  (d+1)^2(d+2)}{d+3}\frac{r_{ps}r_s^{d+1}}{r_{-1/\varepsilon}^2}
\int \phi^2|_{r=r_{-1/\varepsilon}}\,d\sigma\,dt.
\end{multline}
We now use \eqref{rewritten} and the fact that $\frac{13}{12}\le
\frac{d+5}{4}$ for $d\ge 1$ to account for the last term in
\eqref{l_f}.  We then see that
\begin{multline}\label{final}
\int_{\{r\in [r_s,r_{-1/\varepsilon}]\}}
l(f)\phi^2\,r^{d+2}\,dr\,d\omega\,dt 
\ge \int_{\{r\in [r_s,r_{-1/\varepsilon}]\}} l(g)\phi^2\,r^{d+2}\,dr\,d\omega\,dt
\\
+  \frac{1}{48}\int_{\{r\in
    [r_s,r_{-1/\varepsilon}]\}}
\frac{(d+1)^3(d+2)}{d+3}\frac{r_{ps}r_s^{d+1}}{r^4}\frac{1}{r^{d+1}-r_s^{d+1}}a'''(h(r))\phi^2
\,r^{d+2}\,dr\,d\omega\,dt
\\-\frac{13}{18}\int_{\{r\in[r_s,r_{-1/\varepsilon}]\}} 
\Bigl(1-\frac{r_s^{d+1}}{r^{d+1}}\Bigr)^2f'(r)(\partial_r \phi)^2\,r^{d+2}\,dr\,d\omega\,dt
\\-\frac{13}{12}\frac{\delta \varepsilon
  (d+1)^2(d+2)}{d+3}\frac{r_{ps}r_s^{d+1}}{r_{-1/\varepsilon}^2}
\int \phi^2|_{r=r_{-1/\varepsilon}}\,d\sigma\,dt.
\end{multline}
Here we have also used \eqref{df}.  As $l(g)$ remains bounded in the
relevant region and as the coefficient in the integrand of the second term in the right
of \eqref{final} is $\gtrsim \varepsilon^2 e^{1/\varepsilon}\gg 1$ on the
said region, the first term in the right can be controlled by a
fraction of the
second provided $\varepsilon$ is sufficiently small.  This finally
yields
\begin{multline}\label{final2}
\int_{\{r\in [r_s,r_{-1/\varepsilon}]\}}
l(f)\phi^2\,r^{d+2}\,dr\,d\omega\,dt 
\\\ge 
  \frac{1}{96}\int_{\{r\in
    [r_s,r_{-1/\varepsilon}]\}}
\frac{(d+1)^3(d+2)}{d+3}\frac{r_{ps}r_s^{d+1}}{r^4}\frac{1}{r^{d+1}-r_s^{d+1}}a'''(h(r))\phi^2
\,r^{d+2}\,dr\,d\omega\,dt
\\-\frac{13}{18}\int_{\{r\in[r_s,r_{-1/\varepsilon}]\}} 
\Bigl(1-\frac{r_s^{d+1}}{r^{d+1}}\Bigr)^2f'(r)(\partial_r \phi)^2\,r^{d+2}\,dr\,d\omega\,dt
\\-\frac{13}{12}\frac{\delta \varepsilon
  (d+1)^2(d+2)}{d+3}\frac{r_{ps}r_s^{d+1}}{r_{-1/\varepsilon}^2}
\int \phi^2|_{r=r_{-1/\varepsilon}}\,d\sigma\,dt.
\end{multline}
The second term on the right can be bootstrapped into the positive
contribution provided by the first term in the right of \eqref{base}.
The remaining boundary term at $r_{-1/\varepsilon}$ will be controlled
at the end of this section using pieces from the subsequent case.

\noindent{\bf\em Case 2:}
$r_{-1/\varepsilon}
  \le r\le r_{ps}$.

For $r$ in this range, we simply have
\[f(r)=\frac{r^{d+2}-r_{ps}^{d+2}}{r^{d+2}} +
\frac{d+2}{d+3}\frac{r_{ps}r_s^{d+1}}{r^{d+2}}\ln\Bigl(\frac{r^{d+1}-r_s^{d+1}}{\frac{d+1}{2}r_s^{d+1}}\Bigr)\]
which is negative, as is desired in order to guarantee a positive
contribution from the $|\ang \phi|^2$ term of \eqref{base}.  Moreover,
we have
\[f'(r) =
\frac{(d+2)r_{ps}^{d+2}}{r^{d+3}}+\frac{d+2}{d+3}\frac{r_{ps}r_s^{d+1}}{r^2}\frac{d+1}{r^{d+1}-r_s^{d+1}}
-\frac{(d+2)^2}{d+3}\frac{r_s^{d+1}r_{ps}}{r^{d+3}}\ln\Bigl(\frac{r^{d+1}-r_s^{d+1}}{\frac{d+1}{2}r_s^{d+1}}\Bigr)
\]
whose every term is positive for $r_s<r\le r_{ps}$.

It only remains to examine $l(f(r))$ for this region.
Here, we first note that
\begin{equation}\label{lg}l(g(r))=\frac{d+2}{4r^{2d+5}}\Bigl(dr^{2d+2}+(d+3)r_s^{d+1}r^{d+1}-(d+2)^2r_s^{2d+2}\Bigr)
\end{equation}
and
\begin{equation}\label{lh}l\Bigl(\frac{d+2}{d+3}\frac{r_{ps}r_s^{d+1}}{r^{d+2}}h(r)\Bigr)=
-
\frac{(d+2)(d+1)}{4}\frac{r_{ps}r_s^{d+1}}{r^{2d+6}}\bigl(2r^{d+1}-(d+3)r_s^{d+1}\bigr).\end{equation}
Since \eqref{lh} is nonnegative for $r\le r_{ps}$, we have
\[l\Bigl(\frac{d+2}{d+3}\frac{r_{ps}r_s^{d+1}}{r^{d+2}}h(r)\Bigr)\ge
-\frac{(d+2)(d+1)}{4}\frac{r_s^{d+1}}{r^{2d+5}}\bigl(2r^{d+1}-(d+3)r_s^{d+1}\bigr)\]
for $r\le r_{ps}$.
Summing this with \eqref{lg}, we have
\[l(f(r))\ge
\frac{d+2}{4r^{2d+5}}\bigl(r^{d+1}-r_s^{d+1}\bigr)\bigl(dr^{d+1}+r_s^{d+1}\bigr),\]
which is clearly nonnegative for $r\ge r_s$.

\noindent{\bf\em Case 3:} $r_{ps}\le r\le r_\alpha$.

This region corresponds precisely to $h(r)\in [0,\alpha]$.  We, thus,
see that
\[f(r)=\frac{r^{d+2}-r_{ps}^{d+2}}{r^{d+2}}+\frac{d+2}{d+3}\frac{r_{ps}r_s^{d+1}}{r^{d+2}}\frac{h(r)}{15\alpha^4}(3
h(r)^4 -10 h(r)^2\alpha^2 + 15\alpha^4)\]
which is easily seen to be positive.  Moreover,
\[f'(r)=\frac{(d+2)r_{ps}^{d+2}}{r^{d+3}} -
\frac{(d+2)^2r_{ps}r_s^{d+1}}{(d+3)r^{d+3}}a(h(r))
+ \frac{d+2}{d+3}\frac{(d+1)r_{ps}r_s^{d+1}}{r^2(r^{d+1}-r_s^{d+1})}\frac{(h(r)^2-\alpha^2)^2}{\alpha^4}.\]
As $a(h(r))$ takes on a maximum value of $\frac{8\alpha}{15}$ and as
$\frac{(d+2)(d+3)}{2} - \frac{(d+2)^2}{d+3}\frac{8\alpha}{15} > 0$ for
$d\ge 1$ and $\alpha<5$, the sum of the first two terms is positive.
And as the last term is clearly positive, we see that $f'(r)>0$ as
desired.

It remains to verify that $l(f)$ is positive.  We begin by calculating
\begin{multline*}
  l(f)=\frac{d+2}{4r^{2d+5}}\Bigl(dr^{2d+2}+(d+3)r_s^{d+1}r^{d+1}-(d+2)^2
  r_s^{2d+2}\Bigr)
\\+\frac{(d+2)(d+1)r_s^{d+1}r_{ps}}{4(d+3)r^{2d+6}(r^{d+1}-r_s^{d+1})}
\Bigl(-2(d+3)(r^{d+1}-r_{ps}^{d+1})(r^{d+1}-r_s^{d+1})a'(h(r))
\\+(d+1)(d+5)r^{d+1}(r^{d+1}-r_s^{d+1})a''(h(r))
\\-(d+1)^2 r^{2d+2}a'''(h(r))\Bigr).
\end{multline*}
Setting
\begin{align*}
  p(r)&=r(dr^{2d+2}+(d+3)r_s^{d+1}r^{d+1}-(d+2)^2
  r_s^{2d+2})\\
n_1(r)&=
-r_{ps}r_s^{d+1}(d+1)(2r^{d+1}-r_s^{d+1}(d+3))\frac{(h(r)^2-\alpha^2)^2}{\alpha^4}\\
n_2(r)&=r_{ps}r_s^{d+1}\frac{(d+1)^2(d+5)}{d+3}r^{d+1}\frac{4h(r)(h(r)^2-\alpha^2)}{\alpha^4}\\
n_3(r)&=r_{ps}r_s^{d+1}\frac{(d+1)^3}{d+3}\frac{r^{2d+2}}{r^{d+1}-r_s^{d+1}}\cdot 4\frac{\alpha^2-3h(r)^2}{\alpha^4},
\end{align*}
it remains to show that
\[p(r)+n_1(r)+n_2(r)+n_3(r)>0.\]
The dominant term is $p(r)$, and we shall show
\begin{align}
  \frac{1}{3}p(r)+n_1(r)>0,\label{n1}\\
\frac{1}{2}p(r)+n_2(r)\ge 0,\label{n2}\\
\frac{1}{6}p(r)+n_3(r)\ge 0.\label{n3}
\end{align}

\noindent{\em Proof of \eqref{n1}:} Using that we are in the regime
$r\ge r_{ps}$ and that $(h(r)^2-\alpha^2)^2$ is maximized when
$h(r)=0$, we have
\begin{align*}
  \frac{1}{3}p(r)+n_1(r)&\ge
  \frac{1}{3}r_{ps}(dr^{2d+2}+(d+3)r_s^{d+1}r^{d+1}-(d+2)^2r_s^{2d+2})
\\&\qquad\qquad\qquad\qquad -
(d+1)r_{ps}r_s^{d+1}(2r^{d+1}-(d+3)r_s^{d+1})\\
&=\frac{1}{3}r_{ps}\Bigl(dr^{2d+2}-(5d+3)r_s^{d+1}r^{d+1}+(2d^2+8d+5)r_s^{2d+2}\Bigr)\\
&=\frac{1}{3}r_{ps}\Bigl[d\Bigl(r^{d+1}-\frac{5d+3}{2d}r_s^{d+1}\Bigr)^2
  + \frac{1}{4d}(d+1)^2(8d-9)r_s^{2d+2}\Bigr].
\end{align*}
The last quantity is clearly positive, as desired, for $d>1$.

For the case $d=1$, 
\begin{align*}
  \frac{1}{3}p(r)+n_1(r)&\ge
  \frac{1}{3}r(r^4+4r_s^2r^2-9r_s^4)-2\sqrt{2}r_s^3(2r^2-4r_s^2)\\
&=\frac{1}{3}\Bigl(3\sqrt{2}r_s^5 - 13
  r_s^4(r-\sqrt{2}r_s)+20\sqrt{2}r_s^3(r-\sqrt{2}r_s)^2 + 24 r_s^2
  (r-\sqrt{2}r_s)^3
\\&\qquad\qquad\qquad\qquad+5\sqrt{2}r_s(r-\sqrt{2}r_s)^4 +
(r-\sqrt{2}r_s)^5\Bigr)\\
&\ge \frac{1}{3}\Bigl(3\sqrt{2}r_s^5 - 13
  r_s^4(r-\sqrt{2}r_s)+20\sqrt{2}r_s^3(r-\sqrt{2}r_s)^2\Bigr)
\end{align*}
which is an everywhere positive quadratic.

\noindent{\em Proof of \eqref{n2}:}
Here, again, we use that we are studying the regime that $r\ge r_{ps}$
and that $h^2-\alpha^2$ is minimized when $h(r)=0$.  It is thus
obtained that
\begin{equation}\label{n2_1}\begin{split}
  \frac{1}{2}p(r)+n_2(r)&\ge
  \frac{1}{2}r_{ps}(dr^{2d+2}+(d+3)r_s^{d+1}r^{d+1}-(d+2)^2
  r_s^{2d+2})\\
&\qquad\qquad\qquad\qquad
  -\frac{4r_{ps}r_s^{d+1}}{\alpha^2}\frac{(d+1)^2(d+5)}{d+3}r^{d+1}h(r)\\
&=\frac{1}{2}r_{ps}\Bigl(d(r^{d+1}-r_s^{d+1})^2+3(d+1)r_s^{d+1}(r^{d+1}-r_s^{d+1})-(d+1)^2
  r_s^{2d+2}\\
&\qquad\qquad\qquad\qquad
  -\frac{8r_s^{d+1}}{\alpha^2}\frac{(d+1)^2(d+5)}{d+3}r^{d+1}h(r)\Bigr).
\end{split}
\end{equation}
Here, we make the change of variables $x=h(r)$.  Thus,
$r^{d+1}-r_s^{d+1}=\frac{d+1}{2}r_s^{d+1}e^x$.  The right side of
\eqref{n2_1} can be rewritten as
\[\frac{r_{ps}}{2}r_s^{2d+2}(d+1)^2\Bigl[\frac{d}{4} e^{2x} +
  \frac{3}{2}e^x - 1 - \frac{4}{\alpha^2}\frac{(d+5)(d+1)}{d+3}xe^x
  - \frac{8}{\alpha^2}\frac{d+5}{d+3}x\Bigr].\]
Setting 
\[q(x)=\frac{d}{4} e^{2x} +
  \frac{3}{2}e^x - 1 - \frac{4}{\alpha^2}\frac{(d+5)(d+1)}{d+3}xe^x
  - \frac{8}{\alpha^2}\frac{d+5}{d+3}x\]
and noticing that $q(0)>0$, it will suffice to show that $q'(x)\ge 0$
for $x$ between $0$ and $\alpha$.  We compute
\[q'(x)=\frac{1}{2}e^x(3+de^x)-\frac{4}{\alpha^2}\frac{d+5}{d+3}\Bigl[2+(1+d)e^x(1+x)\Bigr].\]
As $\frac{d+5}{d+3}\le \frac{3}{2}$ for $d\ge 1$ and as $1+x\le e^x$,
it follows that
\[q'(x)\ge \frac{1}{2}e^x(3+de^x)-\frac{6}{25}\Bigl(2e^x +
(d+1)e^{2x}\Bigr)>0,\quad d\ge 1,\,\alpha=5.\]
The latter inequality follows from the fact that $\frac{d}{2}\ge
\frac{6}{25}(d+1)$ provided $d\ge 1$.  Since the above inequality
holds for $\alpha = 5$, by continuity, we have that it also
holds for $\alpha=5-\delta_0$ for some $\delta_0>0$, which completes the proof.

%
%
%
%
%
%
%

\noindent{\em Proof of \eqref{n3}:}
Using that we are examining the region $r\ge r_{ps}$ and rewriting
$l(g)$ as in the previous case, we have
\begin{multline*}
  \frac{1}{6}p(r)+n_3(r)\ge\frac{1}{6}r_{ps}\Bigl[d(r^{d+1}-r_s^{d+1})^2+3(d+1)r_s^{d+1}(r^{d+1}-r_s^{d+1})
  - (d+1)^2r_s^{2d+2}\\
 +\frac{24
    r_s^{d+1}}{\alpha^2}\frac{(d+1)^3}{d+3}\Bigl(1-\frac{3}{\alpha^2}h(r)^2\Bigr)\Bigl(
  (r^{d+1}-r_s^{d+1}) + 2 r_s^{d+1} +
  \frac{r_s^{2d+2}}{r^{d+1}-r_s^{d+1}}\Bigr)\Bigr].
\end{multline*}
Using that
\begin{align*}\frac{24r_s^{d+1}}{\alpha^2}\frac{(d+1)^3}{d+3}
\Bigl(1-\frac{3}{\alpha^2}(h(r))^2\Bigr)\frac{r_s^{2d+2}}{r^{d+1}-r_s^{d+1}}
&\ge
-\frac{72r_s^{d+1}}{\alpha^4}\frac{(d+1)^3}{d+3}(h(r))^2\frac{r_s^{2d+2}}{r^{d+1}-r_s^{d+1}}\\
&\ge -\frac{144}{\alpha^4} r_s^{2d+2}\frac{(d+1)^2}{d+3}(h(r))^2\end{align*}
when $r\ge r_{ps}$, we obtain
\begin{multline*}
  \frac{1}{6}p(r)+n_3(r)\ge \frac{1}{6}r_{ps}\Bigl[d(r^{d+1}-r_s^{d+1})^2+3(d+1)r_s^{d+1}(r^{d+1}-r_s^{d+1})
  - (d+1)^2r_s^{2d+2}\\
 +\frac{24
    r_s^{d+1}}{\alpha^2}\frac{(d+1)^3}{d+3}\Bigl(1-\frac{3}{\alpha^2}h(r)^2\Bigr)\Bigl(
  (r^{d+1}-r_s^{d+1}) + 2 r_s^{d+1}\Bigr) \\
  -\frac{144}{\alpha^4}r_s^{2d+2}\frac{(d+1)^2}{d+3}(h(r))^2\Bigr].
\end{multline*}
Proceeding as above with the change of variables $x=h(r)$, this is
\[= \frac{(d+1)^2}{6}r_{ps}r_s^{2d+2}\Bigl[\frac{d}{4}e^{2x}+\frac{3}{2}e^x
  - 1+
 \frac{24}{\alpha^2}\frac{d+1}{d+3}\Bigl(1-\frac{3}{\alpha^2}x^2\Bigr)\Bigl(
  \frac{d+1}{2}e^x + 2\Bigr) 
  -\frac{144}{\alpha^4}\frac{1}{d+3}x^2\Bigr].\]

Setting 
\[s(x)=\frac{d}{4}e^{2x}+\frac{3}{2}e^x
  - 1+
 \frac{24}{\alpha^2}\frac{d+1}{d+3}\Bigl(1-\frac{3}{\alpha^2}x^2\Bigr)\Bigl(
  \frac{d+1}{2}e^x + 2\Bigr) 
  -\frac{144}{\alpha^4}\frac{1}{d+3}x^2,\]
we first note that $s(0)>0$.  For $x\le 5$, we furthermore have
\begin{align*}
  s'(x)&=\frac{1}{2\alpha^4(d+3)}\Bigl[24\alpha^2 (d+1)^2 e^x +
    \alpha^4 (d+3)e^x (3+de^x) - 72 x(8(d+2) + (d+1)^2e^x(x+2))\Bigr]\\
&\ge \frac{1}{2\alpha^4 (d+3)}\Bigl[24\alpha^2 (d+1)^2 e^x + \alpha^4
      (d+3)e^x (3+de^x) - 72 e^x (8(d+2)+7(d+1)^2e^x)\Bigr].
\end{align*}
For $\alpha=5$, this is
\[\frac{1}{1250(d+3)}e^x\Bigl(5073-504 e^x +
51d(49+17e^x)+d^2(600+121e^x)\Bigr),\]
which is easily seen to be positive for $d\ge 1$.  By continuity,
positivity also follows for $\alpha=5-\delta_0$ provided $\delta_0$ is
sufficiently small.

\noindent{\bf\em Case 4:} $r\ge r_\alpha$.

In this regime, 
\[f(r)=\frac{r^{d+2}-r_{ps}^{d+2}}{r^{d+2}}
+\frac{8\alpha}{15}\frac{d+2}{d+3}\frac{r_{ps}r_s^{d+1}}{r^{d+2}},\]
which is clearly positive.  Moreover,
\[f'(r)=\frac{(d+2)r_s^{d+1}r_{ps}}{r^{d+3}}\Bigl(\frac{d+3}{2}-\frac{8\alpha}{15}\frac{d+2}{d+3}\Bigr)\]
which is also positive since $\alpha < 5 \le
\frac{15}{16}\frac{(d+3)^2}{d+2}$ for $d\ge 1$.  Finally, we notice
that $l(f(r))=l(g(r))$ when in this case.  Thus, as in the proof of
\eqref{n2}, we have
\[l(f(r))=r_{ps}\Bigl(d(r^{d+1}-r_s^{d+1})^2+3(d+1)r_s^{d+1}(r^{d+1}-r_s^{d+1})-(d+1)^2
  r_s^{2d+2}\Bigr).\]
As $r^{d+1}-r_s^{d+1}\ge \frac{d+1}{2}r_s^{d+1}$ for $r\ge r_{ps}$, we
see that $l(f(r))>0$ as desired.

\noindent{\bf\em Boundary term at $r_{-1/\varepsilon}$:} 

In order to finish showing that the right side of \eqref{base} is
nonnegative, it remains to examine the $r_{-1/\varepsilon}$ boundary
term in \eqref{base} as well as the subsequent contribution from
\eqref{final2}.  Here, we simply utilize the Fundamental Theorem of
Calculus to control these terms via the positive contributions of the
first and third term in the right of \eqref{base} in the range
$[r_{-1/\varepsilon},r_{ps}]$.  The scaling parameter $\delta$ insures
the necessary smallness.

Fix a smooth cutoff $\beta$ which is identity for, say, $r\le r_{-1}$
and which vanishes for $r\ge r_{ps}$.  Then, for $r\le r_{-1}$, we
have
\[\phi(r)=-\int_r^{r_{ps}} \partial_s (\beta \phi)\,ds.\]
Using the Schwarz inequality, this yields
\[\phi^2(r)\lesssim \int_r^{r_{ps}}|\beta'| \phi^2\,ds -h(r)
\int_r^{r_{ps}}(s^{d+1}-r_s^{d+1}) \beta (\partial_r\phi)^2\,ds.\]
Applying this at $r_{-1/\varepsilon}$ yields
\[\varepsilon \phi^2(r_{-1/\varepsilon})\lesssim
\int_{r_{-1/\varepsilon}}^{r_{ps}}\nabla^\alpha\tilde{P}_\alpha[\phi,X]
r^{d+2}\,dr.\]
Multiplying both sides by $\delta$ and integrating over $[0,T]\times
\S^{d+2}$, we see that these boundary terms can be bootstrapped into
the contributions of {\em Case 2}.

\newsection{A Hardy inequality and the time boundary terms}
In the previous section, we constructed a multiplier so that the right
side of \eqref{base} provides a positive contribution.  By inspection,
the coefficients are easily seen to correspond to those in
\eqref{LEnorm}.  What remains is to control the left side of
\eqref{base} in terms of the initial energy.  For the first term, this
is straightforward.  For the second term in the left side of
\eqref{base}, a Hardy-type inequality is employed, which shall be
proved below.

For the first term in \eqref{base}, we need only apply the Schwarz
inequality to see that
\[\int
f(r)\partial_t\phi(t,\cd)\partial_r\phi(t,\cd)\,r^{d+2}\,dr\,d\omega
\lesssim E[\phi](t).\]
And thus, by conservation of energy, these terms are controlled by
$E[\phi](0)$ as desired.

For the second term in \eqref{base}, we again apply the Schwarz
inequality.  It remains to show that
\begin{equation}\label{time}\int
\Bigl[\frac{1}{r^{d+2}}\partial_r(f(r)r^{d+2})\Bigr]^2\Bigl(1-\frac{r_s^{d+1}}{r^{d+1}}\Bigr)
\phi^2(t,\cd)\,r^{d+2}\,dr\,d\omega \lesssim E[\phi](t)\end{equation}
as a subsequent application of conservation of energy will complete
the proof.

In order to show \eqref{time}, we shall prove a Hardy-type inequality
which is in the spirit of that which appears in \cite{DR}.  Indeed, we
notice that the coefficient in the integrand in the left side of
\eqref{time} is $O((\log(r-r_s))^{-2}(r-r_s)^{-1})$ as $r\to r_s$ and
is $O(r^{-2})$ as $r\to \infty$.  Thus, it will suffice to show that
\begin{equation}\label{hardy}
  \int_{r_s}^\infty
  \frac{1}{r^2}\frac{1}{\Bigl(1-\log\Bigl(\frac{r-r_s}{r}\Bigr)\Bigr)^2\Bigl(\frac{r-r_s}{r}\Bigr)}
  \phi^2\, r^{d+2}\,dr \lesssim \int_{r_s}^\infty
  \Bigl(\frac{r-r_s}{r}\Bigr) (\partial_r\phi)^2\,r^{d+2}\,dr.
\end{equation}

To this end, we set 
\[\rho(r)=\int_{r_s}^r
\frac{x^d}{\Bigl(1-\log\Bigl(\frac{x-r_s}{x}\Bigr)\Bigr)^2\Bigl(\frac{x-r_s}{x}\Bigr)}\,dx.\]
Notice that $\rho(r)\sim r^{d+1}$ as $r\to \infty$ and $\rho(r)\sim
\Bigl[1-\log\Bigl(\frac{r-r_s}{r}\Bigr)\Bigr]^{-1}$ as $r \to r_s$.

Writing the left side of \eqref{hardy} as $\int \rho'(r)\phi^2\,dr$,
integrating by parts, and applying the Schwarz inequality, we have
\begin{align*}
  \int \rho'(r)\phi^2\,dr &= -2\int \rho(r)\phi \partial_r\phi\,dr\\
&\lesssim \Bigl(\int
  \frac{(\rho(r))^2}{\rho'(r)}(\partial_r\phi)^2\,dr\Bigr)^{1/2}
  \Bigl(\int \rho'(r)\phi^2\,dr\Bigr)^{1/2}.
\end{align*}
This completes the proof of \eqref{hardy} as
$\frac{(\rho(r))^2}{\rho'(r)}\sim r^{d+2}$ as $r\to \infty$ and
$\frac{(\rho(r))^2}{\rho'(r)}\sim (r-r_s)$ as $r\to r_s$.

\bigskip

\end{document}